\documentclass[9pt]{amsart}%
\usepackage{amssymb, palatino}
\usepackage{graphicx}
\usepackage{graphicx}%
 \numberwithin{equation}{section}

\newtheorem{dfn}{Definition}[section]
\newtheorem{thm}{Theorem}[section]
\newtheorem{lem}[dfn]{Lemma}
\newtheorem{prop}[dfn]{Proposition}
\newtheorem{rem}[dfn]{Remark}

\newenvironment{preuve}{{\em \bf Proof:}}{\hfill $\square$}
\newenvironment{Ackonowledgement}{{\em  Ackonowledgement:}}

\title{A prescribed Gauss-Kronecker curvature Problem on the product of unit spheres}%
\author{  Wang   Zhizhang }
\date{\today}

\begin{document}

\maketitle

\section{Introduction} \noindent
Prescribed Gauss-Kronecker curvature problems are widely studied in
the literature. Famous among them is the Minkowski problem. It was
studied by H.
 Minkowski, A.D. Alexandrov, H. Lewy, A.V. Pogorelov, L. Nirenberg and at last solved
by S.Y. Cheng and S.T. Yau [CY]. After that, V.I.Oliker [O]
researched the arbitrary hypersurface with prescribed Gauss
curvature in Euclidean space. On the other hand,  L.A. Caffarelli,
L. Nirenberg, and  J. Spruck studied the boundary-value problem of
prescribed Weingarten curvature of graphs over some Euclidean domain
in [CNS2], [CNS3], [CNS4]. Then B. Guan and J. Spruke [GS] studied
the boundary-value problem in the case of hypersurfaces that can be
represented as a radial graph over some domain on some unit sphere.
But on the product unit spheres, the similar problem has not been
studied systematically. The present paper tries to ask and partly
solve a problem of this kind.

Let $S^{m}$ $\subset$ $\mathbb{R}^{m+1}$, $S^{n}$ $\subset$
$\mathbb{R}^{n+1}$, and $S^{m+n+1}$ $\subset$
$\mathbb{R}^{m+1}\oplus\mathbb{R}^{n+1}$ are three unit spheres.
$\vec{\gamma}$, $\vec{\rho}$ are position vectors of $S^m$, $S^n$
respectively, and $u$ is a smooth function defined on $S^m\times
S^n$. Consider a hypersurface $M$ $\subset$ $S^{m+n+1}$ defined by a
natural embedding $\vec{X}$
\begin{eqnarray}
\vec{X}: \ \ \ \ \ S^m\times S^n & \longrightarrow & \ \ \ \ \ \ \ \ \ \ \ \ \ \ \ \ \ \ \ \   S^{m+n+1}\\
\left( \vec{\gamma},\vec{\rho} \right) \ \ \ &
 \longmapsto &
 \frac{1}{(1+e^{-2u})^{1/2}}\vec{\gamma}+\frac{1}{(1+e^{2u})^{1/2}}\vec{\rho}\nonumber.
\end{eqnarray}
 This map
firstly appears in [H], but that paper only discusses the prescribed
mean curvature problem. The fact that the map is an embedding will
be proved in Section 2. Now we state explicitly the main problem:
for a given positive smooth function $K$ defined on $S^m\times S^n$,
can we find a closed strictly convex hypersurface in $S^{m+n+1}$
which is described by (1.1), and whose Gauss-Kronecker curvature is
$K$ ? We will show that there is no global solution to this problem.
So we have to restrict this problem to a subdomain of $S^m\times
S^n$. We solve this problem for some special domains defined as
follows:
\begin{dfn}(PHC-domains) Assume $m\neq n$. For $m>n$, a domain $\Omega$  $\subset$  $S^m\times S^n$
 is called a PHC-domain if it satisfies

\noindent (i) $\Omega$ is a \textbf{product} domain of the form
 $\Omega=\Omega_x\times S^n$ with $\Omega_x$ $\subset$  $S^m$;

\noindent (ii) $\Omega_x$  is  contained in some
\textbf{hemisphere};

\noindent (iii) $\Omega_x$  is  a strictly infinitesimally
\textbf{convex} domain with smooth boundary.

\noindent For $m<n$, we give a similar definition by changing the
position of $m$ and $n$: let $\Omega=S^m\times \Omega_y$, $\Omega_y$
$\subset$  $S^n$, and replace $\Omega_x$ with $\Omega_y$ in
(ii),(iii) above.
\end{dfn}
We know strictly geodesically convex is equivalent to strictly
locally convex, and they can be induced by strictly infinitesimally
convex. For details see [S]. Our main result is the
\begin{thm}
Assume that $\Omega$ $\subset$ $S^m\times S^n$ is a PHC-domain. For
a given smooth positive function $K$ defined on $\bar{\Omega}$,
there is an embedding $\vec{X}$ given by (1.1) on $\bar{\Omega}$,
giving a closed strictly convex smooth hypersurface in $S^{m+n+1}$,
whose Gauss-Kronecker curvature is $K$.
\end{thm}

The prescribed Gauss curvature problems always relate to some
Monge-Amp\`{e}re type equation.  Assume $m$ $>$ $n$. Consider the
Dirichlet problem on a PHC-domain $\Omega$
\begin{eqnarray}
\left\{\begin{array}{ll}\det M(u)\ \ =\ \ K(f(u,|\nabla_xu|^2,|\nabla_yu|^2))^{\frac{m+n+2}{2}}& \text{in}\ \ \ \  \Omega\\
\ \ \ \ \ \ \ \ \ \ \ \ \ \ u\ \ =\ \ \psi & \text{on}  \ \ \
\partial\Omega
\end{array}\right.,
\end{eqnarray}
where $K$, $\psi$ $\in$ $C^{\infty}(\bar{\Omega})$, $K$ $>$ $0$ on
$\bar{\Omega}$, $M(u)$ is defined in (2.9) and the rest of the
notation defined in (2.1),(2.2). This is the equation associated to
our main problem. Although by the effort of much people,
Monge-Amp\`{e}re equations are well understood now, in our case
there are some new difficulties.

The framework to obtain a solution is the  classical continuity
method (see [N]). So we need to give the openness part and the
closedness part. For the openness part, the condition of the
uniqueness of linearized equation relies on the smallness of
boundary-values. Hence we consider Dirichlet problem (1.2) with
sufficiently small boundary-values. For the closedness part, the
openness also gives the comparison Lemma 3.1, leading to the
construction of the subsolution using suitable boundary-values. As
in [CNS1], the subsolution gives the initial solution, the $C^0$
estimate and the $C^1$ estimate on the boundary. The interior $C^1$
estimate is needed, since the manifold we consider here is a product
of a domain and a unit sphere, and the sphere has no boundary. We
choose a natural function (3.22) and estimate it at its maximum
value point to obtain the bound. For the interior $C^2$ estimate,
the difficulty is that we can not diagonalize the three matrices
$(M(u)_{AB})_{(m+n)\times(m+n)}, (M(u)_{ij})_{m\times m}$ and $
(M(u)_{\alpha\beta})_{n\times n} $  at the same time (the
conservation of the indices is stated in the head of Section 2).
Because of this, we need to introduce a term
$\sum_{A}M(u)^{AA}\sum_AM(u)_{AA}$. Inspired by papers of S.T. Yau
[Y] and B. Guan [G], we choose a function (4.11). Then estimating it
at the maximum value point, and computing explicitly the function
$f$ in equation (1.2), we obtain the needed term and the interior
$C^2$ estimate, in which we also generalized the idea of using a
$C^3$ term in the paper [Y]. For the $C^2$ estimate on the boundary,
we have the same difficulty as for the interior estimate and the
difficulty that the manifold is a product manifold. Inspired by
[CNS1], we use a function developed in [G] and the coordinate
functions in Euclidean space to obtain the estimate. Then from
Evans-Krylov theory (see [GT]), we have the $C^{2,\alpha}$ estimate.
At last , differentiate (1.2) and using Schauder theory, we have

\begin{prop}
Let $0$ $<$ $\tau$ $<$ $2$, $K$ $>$ $0$ smooth and $\psi$ $\in$
ABF($\tau$, $K$). Assume that $u$ is the solution of problem (1.2).
Then there is a constant $C_0$ depending on $\psi$, $m$, $n$, $K$,
$\partial\Omega$, $k,\alpha$ such that
\begin{eqnarray}
||u||_{C^{k,\alpha}(\bar{\Omega})}&\leqq& C_0,
\end{eqnarray}
where $k$ is a positive integer and  $0<\alpha<1$.
\end{prop}
Here $\tau$ is defined in (2.1) and ABF($\tau,K$) is defined in
Definition 3.2. Now by the continuity method, we have
\begin{thm}
Let $\Omega\subset S^m\times S^n$ be a PHC-domain, and $m>n$. For a
given smooth function  $K$ $>$ $0$ and $\psi$ $\in$ ABF($\tau$,
$K$), problem (1.2) has a unique convex smooth solution(meaning that
$M(u)$ is positive definite).
\end{thm}
Since the main problem leads to equation (1.2),  for $m>n$ Theorem
1.2 gives the $m>n$ part of Theorem 1.1 . For $m<n$, we take $v=-u$
and change the position of $S^m$ and $S^n$, so  it becomes the
previous case and the map $\vec{X}$ is not changed. Then we obtain
the rest of Theorem 1.1.

The present paper is organized as follows: in Section 2, we compute
out the Gauss-Kroneker curvature of the hypersurface $M$ defined by
map $\vec{X}$, and give the openness part of equation (1.2). Section
3 gives the $C^0$ and $C^1$ estimates of (1.2). And the last two
sections give the $C^2$ estimate in the interior and on the
boundary.

\bigskip

\begin{Ackonowledgement}
The author wishes to thank  Professor JiaXing Hong, Professor
XuanGuo Huang, and Professor Dan Zaffran for their helpful
discussion.
\end{Ackonowledgement}
\section{Equation and Openness} \noindent
We firstly compute the Gauss-Kronecker curvature of the hypersurface
 $M$ defined by (1.1). Let $\{e_1,\cdots,e_m\}$,
 $\{e_{m+1},\cdots,e_{m+n}\}$ be local orthonormal coordinates of
 $S^m$, $S^n$. Throughout our paper, Latin indices ($i,j,\cdots$),
 Greek indices ($\alpha,\beta,\cdots$)
 and capital Latin indices ($A,B,\cdots$) take
 values in the sets $\{1,\cdots,m\}$, $\{m+1,\cdots,m+n\}$ and
 $\{1,\cdots,m+n\}$ respectively. Now we define
\begin{eqnarray}
\tau&=&\frac{2(m-n)}{m+n+2}; \ \ \ \ \ f(r,p,q)\ \ = \ \ e^{\tau
r}\big(1+\frac{p}{1+e^{2r}}+\frac{e^{2r}}{1+e^{2r}}q\big),
\end{eqnarray}
and
\begin{eqnarray}
|\nabla_xu|^2&=&\sum_iu^2_i,\ \ |\nabla_yu|^2\ \ =\ \
\sum_{\alpha}u^2_{\alpha},\ \ |\nabla u|^2\ \ = \ \
|\nabla_xu|^2+|\nabla_yu|^2.
\end{eqnarray}
Obviously,
 \begin{eqnarray}
 \vec{\gamma}_i\ \ =\ \ e_i;\ \ \vec{\gamma}_{\alpha}\ \ =\ \ 0;\ \ \vec{\rho}_i\ \ =\ \ 0;\ \ \vec{\rho}_{\alpha}\ \ =\ \ e_{\alpha}.
 \end{eqnarray}
 Then the tangent vectors of $M$ is
 \begin{eqnarray}
 \vec{X}_A&= &(1+e^{2u})^{-\frac{3}{2}}[e^uu_A\vec{\gamma}+e^u(1+e^{2u})\vec{\gamma}_A+(1+e^{2u})\vec{\rho}_A-e^{2u}u_A\vec{\rho}]
\end{eqnarray}
and the induced metric $g$ is
\begin{eqnarray}
\ \ \ \ \ \ \  \ g_{AB}&=&<\vec{X}_A,\vec{X}_B>\\
&=&\frac{e^{2u}}{(1+e^{2u})^2}[u_Au_B+(1+e^{2u})<\vec{\gamma}_A,\vec{\gamma}_B>+(1+e^{-2u})<\vec{\rho}_A,\vec{\rho}_B>]\nonumber.
\end{eqnarray}
Here $<\cdot, \cdot>$ is the standard inner product of
$\mathbb{R}^{m+n+2}$, and we choose special local coordinates such
that $|\nabla_x u|=u_1$ and $|\nabla_yu|=u_{m+1}$, then
\begin{eqnarray}
\det(g_{AB})&=&\frac{e^{2mu}}{(1+e^{2u})^{m+n}}(1+\frac{|\nabla_xu|^2+e^{2u}|\nabla_yu|^2}{1+e^{2u}}).
\end{eqnarray}
And (2.6) implies that $\vec{X}$ is an embedding. Out of the two
normal unit vectors of $M$ in $S^{m+n+1}$, we choose
\begin{eqnarray}
\vec{n}&=&-\frac{-\vec{\gamma}+e^u\vec{\rho}+\sum_iu_i\vec{\gamma}_i+e^u\sum_{\alpha}u_{\alpha}\vec{\rho}_{\alpha}}
{(1+e^{2u}+|\nabla_xu|^2+e^{2u}|\nabla_yu|^2)^{\frac{1}{2}}},
\end{eqnarray}
and
\begin{eqnarray}
\\
\vec{n}_i&=&
[\ln(1+e^{2u}+|\nabla_xu|^2+e^{2u}|\nabla_yu|^2)^{\frac{1}{2}}]_i\vec{n}
-(1+e^{2u}+|\nabla_xu|^2+
e^{2u}|\nabla_yu|^2)^{-\frac{1}{2}}\nonumber\\
&&\times[\sum_j(u_{ij}-\delta_{ij})\vec{\gamma}_j-u_i\vec{\gamma}
+e^uu_i\vec{\rho}+e^u\sum_{\beta}(u_{i\beta}+u_iu_{\beta})\vec{\rho}_{\beta}]\nonumber,\\
\vec{n}_{\alpha}&=&
[\ln(1+e^{2u}+|\nabla_xu|^2+e^{2u}|\nabla_yu|^2)^{\frac{1}{2}}]_{\alpha}\vec{n}
-(1+e^{2u}+|\nabla_xu|^2+
e^{2u}|\nabla_yu|^2)^{-\frac{1}{2}}\nonumber\\
&&\times[\sum_ju_{j\alpha}\vec{\gamma}_j+e^u\sum_{\beta}(u_{\alpha\beta}+u_{\alpha}u_{\beta}+\delta_{\alpha\beta})\vec{\rho}_{\beta}]\nonumber.
\end{eqnarray}
where we use $\vec{\gamma}_{ij}+\vec{\gamma}\delta_{ij}=0$ and
$\vec{\rho}_{\alpha\beta}+\vec{\rho}\delta_{\alpha\beta}=0$. Now we
denote a symmetric matrix
\begin{eqnarray}
\ \ \ \ \ M(u)_{AB}&=&\left\{\begin{array}{ll}
u_{ij}-u_iu_j-\delta_{ij} & A,B\in \{1,\cdots,m\}\\ u_{i\alpha} &
A\in \{1,\cdots,m\},
B\in\{m+1,\cdots,m+n\}\\
u_{\alpha\beta}+u_{\alpha}u_{\beta}+\delta_{\alpha\beta} &
A,B\in\{m+1,\cdots,m+n\}\end{array}\right..
\end{eqnarray}
Then the second fundamental tensor of $M$ along $\vec{n}$ is
\begin{eqnarray}
h_{AB}&=&-<\vec{n}_A,\vec{X}_B>\\
&=&e^u[(1+e^{2u})(1+e^{2u}+|\nabla_x
u|^2+e^{2u}|\nabla_yu|^2)]^{-\frac{1}{2}}M(u)_{AB}.\nonumber
\end{eqnarray}
The Gauss-Kronecker curvature of $M$ in $S^{m+n+1}$ is
\begin{eqnarray}
K&=&\frac{\det(h_{AB})}{\det(g_{AB})}\ \ =\
 \ \frac{e^{(n-m)u}}{(1+\dfrac{|\nabla_xu|^2+e^{2u}|\nabla_yu|^2}{1+e^{2u}})^{\frac{m+n+2}{2}}}
\det M(u).
\end{eqnarray}
Then rewrite (2.11), and using (2.1),(2.2), we have (1.2).
\begin{prop}
There is no global strictly convex hypersurface described by (1.1)
on the product of unit spheres.
\end{prop}
\begin{preuve} The meaning of strictly convex here is that the
hypersurface $M$ viewed as the submanifold of $S^{m+n+1}$ is
strictly convex. This means that the second fundamental tensor is
positive or negative definite, which implies that
$(u_{ij}-u_iu_j-\delta_{ij})$ $>$ $0$ or
$(u_{\alpha\beta}+u_{\alpha}u_{\beta}+\delta_{\alpha\beta})$ $<$
$0$. In the first case, at the maximum value point of $u$, we have
$(u_{ij})$ $\leqq$ $0$ and $\nabla u=0$. Then $(\delta_{ij}) <
(u_{ij}) \leqq 0$, which is a contradiction. The second case has the
same contradiction at the minimum value point of $u$.
\end{preuve}

\bigskip

\begin{rem}
In the problem (1.2), $u$ attains its maximum value point only on
the boundary of  $\Omega$. We can see this by the same argument as
used in Proposition 2.1.
\end{rem}
\hfill$\square$

Now we discuss the openness part of  problem (1.2). Let $D=\{u\in
C^{2,\alpha}(\bar{\Omega});u|_{\partial\Omega}=\psi, M(u)\text{ is
positive definite} \}$. Consider the map $F$ induced by problem
(1.2)
\begin{eqnarray}
F: D&\longrightarrow&\
\ \ \ \ \ \ \ \ \ \ \ \ \ \ \ \ C^{\alpha}(\bar{\Omega})\\
u  &\longmapsto&
(f(u,|\nabla_xu|^2,|\nabla_yu|^2))^{-\frac{m+n+2}{2}}\det
M(u)\nonumber.
\end{eqnarray}
The linearized operator of $F$ at $u$ along function $v$ is
\begin{eqnarray}
&&DF(u)v\\
&=&f^{-\frac{m+n+2}{2}}\det M(u)[\sum_{A,B}M(u)^{AB}v_{AB}-2\sum_{i,j}M(u)^{ij}u_iv_j+2\sum_{\alpha,\beta}M(u)^{\alpha\beta}u_{\alpha}v_{\beta}]\nonumber\\
&&-\frac{m+n+2}{2}f^{-\frac{m+n+4}{2}}\det
M(u)[f_rv+2f_p\sum_iu_iv_i+2f_q\sum_{\alpha}u_{\alpha}v_{\alpha}]\nonumber,
\end{eqnarray}
where $(M(u)^{AB})$ is the inverse matrix of  $(M(u)_{AB})$. When
$f_r$ $\geqq$ $0$, the linearized problem has a unique solution.
Then by linear elliptic PDE theory, $DF(u)$ is a continuous  linear
bijective map. Now by the implicit function theorem and the openness
of positivity, we obtain the openness. Hence, we only need to find
the condition guaranteeing $f_r$ $\geqq$ $0$. Since
\begin{eqnarray}
f_r&=&e^{\tau
r}\big[\tau+\frac{\tau(1+e^{2r})-2e^{2r}}{(1+e^{2r})^2}p+\frac{\tau
e^{2r}q}{1+e^{2r}}+\frac{2e^{2r}q}{(1+e^{2r})^2}\big],
\end{eqnarray}
at $r=u$, $p=|\nabla_xu|^2$, $q=|\nabla_yu|^2$, we only need
\begin{eqnarray}
\tau(1+e^{2u})-2e^{2u}&\geqq& 0.
\end{eqnarray}
Now by Remark 2.2, if we require
\begin{eqnarray}
u\ \ \leqq\ \  \sup_{\partial \Omega} u\ \ =\ \ \sup_{\partial
\Omega}\psi\ \ \leqq \ \ \frac{1}{2}\ln\frac{\tau/2}{1-\tau/2},
\end{eqnarray}
then (2.15) is satisfied. And (2.16) makes sense if $m$ $>$ $n$ by
(2.1).
\begin{prop}
Assume $m$ $>$ $n$, and $\sup_{\partial \Omega}$ $\psi$ $\leqq$
$\dfrac{1}{2}\ln\dfrac{\tau/2}{1-\tau/2}$, then the openness of
problem (1.2) holds.
\end{prop}
\hfill$\square$
\section{$C^0$ and $C^1$ estimates} \noindent
From this section on, we always assume $m$ $>$ $n$ and use the
Einstein convention: indices appearing twice in an expression, once
as a subscript,once as a superscript implicitly summed over. The
following comparison Lemma maybe known, but the author did not find
an appropriate reference, so it is included here.
\begin{lem}
Let $\Omega$ be a domain on the product of unit spheres. Let
\begin{eqnarray}
G(u)&=&\det M(u)[g(u,|\nabla_x u|^2,|\nabla_yu|^2)]^{-1},
\end{eqnarray}
where $M(u)$ defined in (2.9) is positive definite and $g(r,p,q)$ is
a positive smooth function on $\mathbb{R}^3$. Assume that there is a
constant $r_0$ such that for $r\leqq r_0$, $g_r>0$. For two smooth
functions $u,v$, if $ u,v\leqq r_0$ and $G(u)\leqq G(v)$ then one of
the following holds
\begin{eqnarray}
v-\sup_{\partial\Omega}(v-u)\ \ \leqq \ \ u,\ \ \ \   \text{or}  \ \
\ \  v\ \ \leqq  \ \ u.
\end{eqnarray}
\end{lem}
\begin{preuve}
By $G(u)\leqq G(v)$
\begin{eqnarray}
\frac{\det M(v)}{\det M(u)}&\geqq&\frac{g(v,|\nabla_x
v|^2,|\nabla_yv|^2)}{g(u,|\nabla_x u|^2,|\nabla_yu|^2)}.
\end{eqnarray}
Since $M(u)$ is symmetric and positive definite, we can assume
$M(u)$ $=$ $CC^T$, where $C$ is a non-degenerate matrix. Then
\begin{eqnarray}
\frac{\det M(v)}{\det M(u)}&=&\det[C^{-1}(M(v)-M(u))(C^T)^{-1}+Id]\\
&\leqq&
\big[\frac{m+n+tr(C^{-1}(M(v)-M(u))(C^{-1})^T)}{m+n}\big]^{m+n}\nonumber\\
&=&\big\{1+\frac{1}{m+n}tr[M(u)^{-1}(M(v)-M(u))]\big\}^{m+n}\nonumber,
\end{eqnarray}
where, $tr$ means taking the trace of a matrix. If (3.2) does not
hold, and the function $v-u$ attains its maximum value at some point
in $\Omega$, then at that point, $v-u$ $>$ $0$,
 $\nabla v=\nabla u$, and $[(v-u)_{AB}]$ is
non-positive definite. So
\begin{eqnarray}
tr[(M(u))^{-1}(M(v)-M(u))]& =& M(u)^{AB}(v-u)_{AB}\ \ \leqq \ \ 0,
\end{eqnarray}
And since $u<v\leqq r_0$, $\nabla u=\nabla v$, $g_r>0$, we have
\begin{eqnarray}
g(v,|\nabla_x
v|^2,|\nabla_yv|^2)&=&g(v,|\nabla_xu|^2,|\nabla_yu|^2)\ \ > \ \
g(u,|\nabla_x u|^2,|\nabla_yu|^2).
\end{eqnarray}
Hence by (3.3)-(3.6), we have a contradiction.
\end{preuve}
\bigskip

We introduce a class of function sets,
\begin{dfn}(ABF-sets) Assume $\Omega$ is a PHC-domain, $m>n$
 and $K$ is a smooth positive function defined on
 $\bar{\Omega}$.  We assume $\Omega=\Omega_x\times S^n$.
Let $\psi$ be a smooth function on $\bar{\Omega}_x$. We call $\psi$
an \textbf{admissible boundary function} (simply, ABF) with respect
to $\tau$ and $K$, if it satisfies

\noindent (i) The matrix $M(\psi)$ is positive definite,

\noindent (ii) $\sup_{\partial\Omega}\psi$ $\leqq$
$\dfrac{1}{2}\ln\dfrac{\tau/3}{1-\tau/2}$, and $F(\psi)$ $\geqq$
$K$.

\noindent The set of all ABFs with respect to $\tau$ and $K$ is
denoted by ABF($\tau$, $K$).
\end{dfn}

In fact for $0$ $<$ $\tau$ $<$ $2$, $K$ $>$ $0$ on a PHC-domain,
 the set ABF($\tau$, $K$) is always non-empty. Indeed for a PHC-domain $\Omega$, we can always assume that
$\bar{\Omega}_x$ is contained in the hemisphere $\{ x\in
\mathbb{R}^{m+1};x\in S^m \ \ \ and\ \ \ x_{1}>0\}$. So we can take
a constant $E$ $>$ $0$ sufficiently small, such that $x_{1}-E$ $>$
$0$. Then for any constant $F$ $>$ $0$, let
\begin{eqnarray}
\varphi&=&-\ln F(x_{1}-E).
\end{eqnarray}
The derivative of $\varphi$ is
\begin{eqnarray}
\varphi_i&=&-\frac{(x_{1})_i}{x_{1}-E},\ \ \ \ \varphi_{ij}\ \ = \ \
\frac{x_{1}\delta_{ij}}{x_{1}-E}+\frac{(x_{1})_i(x_{1})_j}{(x_{1}-E)^2},
\end{eqnarray}
where we use $(x_{1})_{ij}$ $+$ $x_{1}\delta_{ij}$ $=$ $0$
(cf.[CY]). Then
\begin{eqnarray}
\varphi_{ij}-\varphi_i\varphi_j-\delta_{ij}&=&\frac{E}{x_{1}-E}\delta_{ij}.
\end{eqnarray}
Then we take $A$ big enough such that
\begin{eqnarray}
\varphi-A &\leqq & \frac{1}{2}\ln\dfrac{\tau/3}{1-\tau/2},
\end{eqnarray}
and
\begin{eqnarray} F(\varphi-A)\ \ \geqq \ \
e^{(m-n)(A-\varphi)}\frac{\det M(\varphi)}{(1+|\nabla
\varphi|^2)^{\frac{(m+n+2)}{2}}} \ \ \geqq \ \ K
\end{eqnarray}
both holds. So let $\psi=\varphi-A$, then $\psi$ satisfies (i)(ii)
of definition (3.2). Then ABF($\tau$, $K$) is always not empty.

Now fix $\tau$, $\tilde{K}$. For any $\tilde{\psi}$ $\in$
ABF($\tau$, $\tilde{K}$) as the boundary function, we use the
continuity method to solve the Dirichlet problem (1.2). For $t$
$\in$ $[0,1]$, take
\begin{eqnarray}
\tilde{K}_t&=&(1-t)F(\tilde{\psi})+t\tilde{K}.
\end{eqnarray}
Assume that $F(u_t)=\tilde{K}_t$  with
$u_t|_{\partial\Omega}=\tilde{\psi}$. Obviously, at $t=0$, we can
take $u_0=\tilde{\psi}$. By Definition 3.2, $F(\tilde{\psi})\geqq
\tilde{K}$, so
\begin{eqnarray}
\tilde{K}&\leqq&\tilde{K}_t\ \ \leqq \ \ F(\tilde{\psi}).
\end{eqnarray}
Let function $g$ equals $f^{(m+n+2)/2}$ in Lemma 3.1. Take $r_0=
\dfrac{1}{2}\ln\dfrac{\tau/3}{1-\tau/2}$. When $r\leqq r_0$, we have
\begin{eqnarray}
\tau(1+e^{2r})-2e^{2r}&\geqq& \frac{\tau}{3}.
\end{eqnarray}
By (2.14) and $0<\tau<2$, for $p,q>0$,
\begin{eqnarray}
\frac{f_r(r,p,q)}{f(r,p,q)}&\geqq &
\dfrac{\tau+\dfrac{\dfrac{\tau}{3}p}{(1+e^{2r})^2}+\dfrac{2e^{2r}q}{(1+e^{2r})^2}}{1+\dfrac{p}{1+e^{2r}}+\dfrac{e^{2r}q}{1+e^{2r}}}
\ \ \geqq \ \ \frac{\tau/3}{1+e^{2r}}\ \ >\ \ 0.
\end{eqnarray}
 Since we have $\tilde{\psi}$
$\in$ ABF($\tau,\tilde{K})$, by the (ii) of Definition 3.2 we have
$u_t\leqq r_0$. Then (3.15) implies $g_r>0$. Moreover by (3.13),
$F(u_t)\leqq F(\tilde{\psi})$. Now we can use Lemma 3.1. By
$u_t|_{\partial\Omega}=\tilde{\psi}$, we have
\begin{eqnarray}
\tilde{\psi}&\leqq& u_t.
\end{eqnarray}
Combining (3.16) and Remark 2.2, we obtain the $C^0$ estimate.

 Now by the definition of derivative, (3.16) and $u_t|_{\partial\Omega}$ $=$ $\tilde{\psi}$, on $\partial\Omega$,
\begin{eqnarray}
\frac{\partial u_t}{\partial n}&\leqq& \frac{\partial
\tilde{\psi}}{\partial n},
\end{eqnarray}
where $\vec{n}$ is the outer normal direction of $\Omega$. If vector
$Y$ is in the tangent space of the submanifold $\partial\Omega$,
then $\nabla_Yu_t= \nabla_Y\tilde{\psi}$. Hence there is a constant
$C_1$ depending on
 $\tilde{\psi}$, $m$, $n$, such that
\begin{eqnarray}
\nabla_Yu_t&\leqq& C_1,
\end{eqnarray}
where $Y$ is a unit vector in the tangent space of $S^m\times S^n$
supported by $\partial \Omega$, and the angle between $Y$ and
$\vec{n}$ is not bigger than $\pi /2$.

For any point $P$ $\in$ $\partial\Omega_x$ $\times$
$S^n$, we know that $P$ $=$ $(P_x,P_y)$, where $P_x$ $\in$ $
\partial\Omega_x$, $P_y$ $\in$ $S^n$.  In $\Omega_x$ we take a geodesic curve
$l_x$ which starts at $P_x$ with direction $-\vec{n}(P_x)$ ( the
inward normal direction of $P_x$ ) and ends at point $Q_x$ $\in$
$\partial\Omega_x$. Denote by $\sigma$ the arc parameter of $l_x$.
Then at every point of $l_x$, we choose a orthonormal frame
$\{e_1,\cdots,e_m\}$ such that $e_m$ is the positive direction of
$l_x$, namely $e_m=\partial/(\partial \sigma)$. Then by the
positivity of the matrix $M(u)$ along geodesic curve $l_x$,
\begin{eqnarray}
(u_t)_{\sigma;\sigma}-(u_t)^2_{\sigma}-1&>&0,
\end{eqnarray}
where $(u_t)_{\sigma,\sigma}$ is the second order normal derivative.
Take integral and by (3.18),
\begin{eqnarray}
-\frac{\partial u_t}{\partial n}(P)&=&(u_t)_{\sigma}(P_x,P_y)\ \  <\
\ (u_t)_{\sigma}(Q_x,Q_y)\ \ \leqq \ \  C_1.
\end{eqnarray}
Combining (3.18) (3.20), we obtain the $C^1$  estimate on the
boundary. Namely, there is a constant $C_2$ depending on
$\tilde{\psi}$, $m$, $n$, $\partial\Omega$, such that on $\partial
\Omega$
\begin{eqnarray}
|\nabla u_t|&\leqq&C_2.
\end{eqnarray}

Now we only need to give the interior $C^1$ estimate. Without loss
of the generality, we only give the estimate for equation (1.2).

Consider a function
\begin{eqnarray}
\phi&=&|\nabla u|^2,
\end{eqnarray}
and let the operator
\begin{eqnarray}
Lv&=&M(u)^{AB}v_{AB}.
\end{eqnarray}
Assume that $\phi$ attains its maximum value at some point $P$ in
$\Omega$. Then at $P$,
\begin{eqnarray}
\phi_A&=&2\sum_Cu_Cu_{CA}=0.
\end{eqnarray}
By the Ricci identity on the product of unit spheres,
\begin{eqnarray}
u_{ABC}-u_{ACB}&=&-u_C\delta_{AB}+u_B\delta_{AC}.
\end{eqnarray}
Using (3.25) and the positivity of $M(u)$,
\begin{eqnarray}
\frac{1}{2}L\phi&=&
\sum_CM(u)^{AB}u_{CA}u_{CB}+\sum_CM(u)^{AB}u_Cu_{CAB}\\
&\geqq&\sum_CM(u)^{AB}u_Cu_{ABC}+\sum_AM(u)^{AA}|\nabla
u|^2-M(u)^{AB}u_Au_B\nonumber\\
&\geqq& \sum_CM(u)^{AB}u_Cu_{ABC}\nonumber.
\end{eqnarray}
Take logarithm of (1.2), and differentiate it. We get
\begin{eqnarray}
M(u)^{AB}u_{ABC}&=&M(u)^{ij}(u_iu_j+\delta_{ij})_C-M(u)^{\alpha\beta}(u_{\alpha}
u_{\beta}+\delta_{\alpha\beta})_C+\frac{K_C}{K}\\
&&+\frac{m+n+2}{f}\big(\frac{f_r}{2}u_C+f_p\sum_iu_iu_{iC}+f_q\sum_{\alpha}u_{\alpha}u_{\alpha
C}\big)\nonumber.
\end{eqnarray}
Now by (3.24) and (3.26),
\begin{eqnarray}
L\phi&\geqq&2\sum_C\frac{K_Cu_C}{K}+\frac{m+n+2}{f}f_r|\nabla u|^2.
\end{eqnarray}
Then by  (3.15) and the $C^0$ estimate, there is a positive constant
$\varepsilon_1$ depending on $\psi$, $m$, $n$, such that
\begin{eqnarray}
\frac{f_r(u,|\nabla_x u|^2,|\nabla_y u|^2)}{f(u,|\nabla_x
u|^2,|\nabla_y u|^2)}\geqq  \varepsilon_1.
\end{eqnarray} Now by
(3.29) and Schwarz inequality, at point $P$
\begin{eqnarray}
0\ \ \geqq \ \ L \phi \ \ \geqq \ \ -C-2\varepsilon_1|\nabla
u|^2+(m+n+2)\varepsilon_1|\nabla u|^2,
\end{eqnarray}
which implies the interior estimate. Here $C$ is a positive constant
depending on $\psi,K,m,n$. Now we have proved
\begin{prop}
Let $0$ $<$ $\tau$ $<$ $2$, $\tilde{K}$ $>$ $0$ and $\tilde{\psi}$
$\in$ ABF($\tau$, $\tilde{K}$). Assume that $u_t$ is the solution of
problem (1.2) in which $K$ is $\tilde{K}_t$ defined by (3.12), and
$\psi$ is $\tilde{\psi}$. Then there is a constant $C_3$ depending
on $\tilde{\psi}$, $m$, $n$, $\tilde{K}$, $\partial\Omega$ such that
\begin{eqnarray}
||u_t||_{C^1(\bar{\Omega})}&\leqq& C_3.
\end{eqnarray}
\end{prop}
\hfill$\square$
\section{ Interior $C^2$ estimate}\noindent
Since we have a lot of positive constants, for simplicity, from this
section on, we write $C$ to represent any constant of minor
important. For a useful constant, we use $C$ or $\hat{C}$ with a
lower index (for example $C_1$, $\hat{C}_1$) to represent it. These
constants always relate to $\psi,K,m,n$ and the $C^1$ norm of $u$,
but we do not refer to this fact everywhere. Without loss of
generality, we only estimate problem (1.2) with $\psi\in$
ABF($\tau,K$), so by (3.14), (3.15) and (3.16), $u\geqq \psi$. Since
$\Omega$ is a PHC-domain, we can assume $\bar{\Omega}_x$ $\subset$
$\{ x\in \mathbb{R}^{m+1};x\in S^m \ \ \ and\ \ \ x_{1}>0\}$. Then
there is a positive constant $\varepsilon_2$ depending on
$\partial{\Omega}$ such that in $\Omega$,
\begin{eqnarray}
x_1&\geqq& \varepsilon_2.
\end{eqnarray}
We let
\begin{eqnarray}
\eta&=& 3-e^{-C_4(u-\psi+1)}-e^{-C_{5}x_1},\\
\zeta&=&3-\eta,
\end{eqnarray}
where $C_4,C_5$ are two positive constants which will be determined
in the following. Since the matrix $M(\psi)$ is positive, we can
assume $M(\psi)\geqq 4\varepsilon_3id$, where $\varepsilon_3$ is a
positive constant depending on $\psi$. Then using (2.9),(3.23) and
Proposition 3.3,
\begin{eqnarray}
&&\ \ L(e^{-C_4(u-\psi+1)})\\
&=&e^{-C_4(u-\psi+1)}\{C_4^2M(u)^{AB}(u-\psi)_A(u-\psi)_B-C_4M(u)^{AB}(M(u)_{AB}-M(\psi)_{AB})\nonumber\\
&&-C_4M(u)^{ij}[u_iu_j-\psi_i\psi_j]-C_4M(u)^{\alpha\beta}[-u_{\alpha}u_{\beta}+\psi_{\alpha}\psi_{\beta}]\}\nonumber\\
&\geqq&
e^{-C_4(u-\psi+1)}\{C_4^2M(u)^{AB}(u-\psi)_A(u-\psi)_B-(m+n)C_4+4\varepsilon_3C_4\sum_AM(u)^{AA}\nonumber\\
&&-CC_4\sum_iM(u)^{ii}+C_4M(u)^{\alpha\beta}(u-\psi)_{\alpha}u_{\beta}+C_4M(u)^{\alpha\beta}\psi_{\alpha}(u-\psi)_{\beta}\}\nonumber.
\end{eqnarray}
Now
\begin{eqnarray}
&&M(u)^{AB}(-\eta)_{AB}+2M(u)^{ij}u_i\eta_{j}-2M(u)^{\alpha\beta}u_{\alpha}\eta_{\beta}\\
&\geqq&
e^{-C_4(u-\psi+1)}\{C_4^2M(u)^{AB}(u-\psi)_A(u-\psi)_B-(m+n)C_4+4\varepsilon_3C_4\sum_AM(u)^{AA}\nonumber\\
&&-CC_4\sum_iM(u)^{ii}+3C_4M(u)^{\alpha\beta}(u-\psi)_{\alpha}u_{\beta}+C_4M(u)^{\alpha\beta}\psi_{\alpha}(u-\psi)_{\beta}\}\nonumber\\
&&+e^{-C_{5}x_1}[C_{5}^2M(u)^{ij}(x_1)_{i}(x_1)_j+2C_{5}M(u)^{ij}u_{i}(x_1)_{j}-C_{5}M(u)^{ij}(x_1)_{ij}]\nonumber\\
&=&e^{-C_4(u-\psi+1)}\{C_4^2M(u)^{AB}(u-\psi)_A(u-\psi)_B-(m+n)C_4+4\varepsilon_3C_4\sum_AM(u)^{AA}\nonumber\\
&&-CC_4\sum_iM(u)^{ii}+3C_4M(u)^{\alpha
B}\psi_{\alpha}(u-\psi)_{B}+C_4M(u)^{A\beta}(u-\psi)_{A}u_{\beta}\nonumber\\
&&-3C_4M(u)^{\alpha
i}\psi_{\alpha}(u-\psi)_i-C_4M(u)^{i\beta}(u-\psi)_iu_{\beta}\}\nonumber\\
&&+e^{-C_{5}x_1}[C_{5}^2M(u)^{ij}(x_1)_{i}(x_1)_j+2C_{5}M(u)^{ij}u_{i}(x_1)_{j}-C_{5}M(u)^{ij}(x_1)_{ij}]\nonumber.
\end{eqnarray}
By Schwarz inequality and because the $C^1$ norms of $\psi$ and $u$
are bounded, we have
\begin{eqnarray}
|3C_4M(u)^{A\beta}(u-\psi)_{A}u_{\beta}|
&\leqq&\frac{C_4^2}{2}M(u)^{AB}(u-\psi)_A(u-\psi)_B+\frac{9}{2}M(u)^{\alpha\beta}u_{\alpha}u_{\beta}\\
&\leqq&\frac{C_4^2}{2}M(u)^{AB}(u-\psi)_A(u-\psi)_B+C\sum_AM(u)^{AA}\nonumber,
\end{eqnarray}
and we also have a similar inequality for the term $C_4M(u)^{\alpha
B}(\psi)_{\alpha}(u-\psi)_{B}$ and the inequality
$|2C_{5}M(u)^{ij}u_{i}(x_1)_{j}|\leqq
C_5^2M(u)^{ij}(x_1)_i(x_1)_j+M(u)^{ij}u_iu_j$. Then
\begin{eqnarray}
&&M(u)^{AB}(-\eta)_{AB}+2M(u)^{ij}u_i\eta_{j}-2M(u)^{\alpha\beta}u_{\alpha}\eta_{\beta}\\
&\geqq&e^{-C_4(u-\psi+1)}\{-\hat{C}_1C_4+(4\varepsilon_3C_4-\hat{C}_2)\sum_AM(u)^{AA}-\hat{C}_{3}C_4\sum_iM(u)^{ii}\nonumber\\
&&-2\hat{C}_{4}C_4\sum_{i,\alpha}|M(u)^{\alpha i}|\}+e^{-C_{5}x_1}[-M(u)^{ij}u_iu_j-C_{5}M(u)^{ij}(x_1)_{ij}]\nonumber.\nonumber,
\end{eqnarray}
where $\hat{C}_1,\cdots,\hat{C}_4$ are four positive constants.
Obviously by the positivity of $M(u)$,
\begin{eqnarray}
\ \ \ \ 2|M(u)^{\alpha i}|&\leqq&2|M(u)^{\alpha
\alpha}M(u)^{ii}|^{1/2}\ \ \leqq\ \
\frac{\varepsilon_3}{m\hat{C}_{4}}M(u)^{\alpha
\alpha}+\frac{m\hat{C}_{4}}{\varepsilon_3}M(u)^{ii}.
\end{eqnarray}
Now by (4.7) and the fact $(x_1)_{ij}+x_1\delta_{ij}=0$,
\begin{eqnarray}
&&M(u)^{AB}(-\eta)_{AB}+2M(u)^{ij}u_i\eta_{j}-2M(u)^{\alpha\beta}u_{\alpha}\eta_{\beta}\\
&\geqq&e^{-C_4(u-\psi+1)}[-\hat{C}_1C_4+(3\varepsilon_3C_4-\hat{C}_2)\sum_AM(u)^{AA}-\hat{C}_{5}C_4\sum_iM(u)^{ii}]\nonumber\\
&&+e^{-C_{5}x_1}(C_{5}x_1-C_3^2)\sum_{i}M(u)^{ii}\nonumber,
\end{eqnarray}
where $\hat{C}_5$ is also a positive constant depending on
$\hat{C}_4, \varepsilon_3,m$ and $n$. Now we take
$C_{5}=2C_3^2/\varepsilon_2$, and
$C_4=\max\{\hat{C}_2/\varepsilon_3,2(\hat{C}_5e^{C_{5}}/ C^2_3)\}$.
Now since $ x_1\leqq 1, u\geqq \psi$ and (4.1), we have
\begin{eqnarray}
M(u)^{AB}(-\eta)_{AB}+2M(u)^{ij}u_i\eta_{j}-2M(u)^{\alpha\beta}u_{\alpha}\eta_{\beta}&\geqq&4\varepsilon_4\sum_AM(u)^{AA}-C,
\end{eqnarray}
where $\varepsilon_4$ $=$
$e^{-C_4(1+\max_{\bar{\Omega}}|u-\psi|)}\hat{C}_2/2$. Now define a
function
\begin{eqnarray}
\phi&=& e^{-C_6\eta}\Delta u+\zeta,
\end{eqnarray}
where $C_6$  is a positive constant which will be determined in the
following, and $\Delta$ is the Laplace operator of $S^m\times S^n$.
This type of function is well known (see [Y]), but we modify it and
use the idea of [B] to handle the extra term
$\sum_{A}M(u)^{AA}\sum_AM(u)_{AA}$ which will appear in the
following. Assume $\phi$ attains its maximum value at point $P$
$\in$ $\Omega$. Then at $P$,
\begin{eqnarray}
\phi_A &=&-C_6\eta_Ae^{-C_6\eta}\Delta u+e^{-C_6\eta}\sum_Cu_{CCA}+\zeta_A\ \ = \ \ 0.
\end{eqnarray}
Then
\begin{eqnarray}
&&L(\phi)\\ &=&e^{-C_6\eta}[C_6^2M(u)^{AB}\eta_A\eta_B\Delta
u-C_6M(u)^{AB}\eta_{AB}\Delta u-2C_6M(u)^{AB}\eta_A\sum_C
u_{CCB}\nonumber\\&&+\sum_C
M(u)^{AB}u_{CCAB}]+M(u)^{AB}\zeta_{AB}\nonumber.
\end{eqnarray}
Use the notation $R^A_{BCD}$ to denote the Riemannian curvature(see
[C] appendix A.6 and the Ricci identities using in the following
also see this book). Firstly by the Ricci identity and (2.9),
\begin{eqnarray}
&&M(u)^{AB}\eta_A u_{CCB}\\&=&
M(u)^{AB}\eta_Au_{BCC}-M(u)^{AB}\eta_Au_ER^E_{BCC}\nonumber\\
&\leqq&M(u)^{AB}\eta_AM(u)_{BCC}+M(u)^{Ai}\eta_A(u_iu_j+\delta_{ij})_j-M(u)^{A\alpha}\eta_A(u_{\alpha}u_{\beta}
+\delta_{\alpha\beta})_{\beta}\nonumber\\&&+ C\sum_A
M(u)^{AA}\nonumber\\
&=& M(u)^{AB}\eta_AM(u)_{BCC}+M(u)^{Ai}\eta_Au_i\sum_j
M(u)_{jj}+M(u)^{Ai}\eta_Au_i\sum_j(u_j^2+1)\nonumber\\
&&-M(u)^{A\alpha}\eta_Au_{\alpha}\sum_{\beta}M(u)_{\beta\beta}+M(u)^{A\alpha}\eta_Au_{\alpha}\sum_{\beta}(u_{\beta}^2+1)
+M(u)^{Ai}\eta_Au_{ij}u_j \nonumber\\
&&-M(u)^{A\alpha}\eta_Au_{\alpha\beta}u_{\beta}
+C\sum_AM(u)^{AA}\nonumber\\
&\leqq&M(u)^{AB}\eta_AM(u)_{BCC}+(|M(u)^{Ai}\eta_Au_i|+|M(u)^{A\alpha}\eta_Au_{\alpha}|)\sum_AM(u)_{AA}\nonumber\\
&&+M(u)^{Ai}\eta_Au_{ij}u_j-M(u)^{A\alpha}\eta_Au_{\alpha\beta}u_{\beta}+C\sum_AM(u)^{AA}\nonumber.
\end{eqnarray}
Since
\begin{eqnarray}
|M(u)^{Ai}\eta_Au_i|&\leqq& \frac{\varepsilon_4}{2}\sum_i M(u)^{ii}+
CM(u)^{AB}\eta_A\eta_B,
\end{eqnarray}
and we have a similar inequality for the term
$|M(u)^{A\alpha}\eta_Au_{\alpha}|$. By Schwarz inequality and (2.9)
\begin{eqnarray}
&&M(u)^{Ai}\eta_Au_{ij}u_j\ \ \leqq \ \
M(u)^{Ai}\eta_A(M(u)_{ij}u_j)+C\sum_AM(u)^{AA}\\
&\leqq&
(M(u)^{AB}\eta_A\eta_B)^{1/2}(M(u)^{ij}M(u)_{ik}u_kM(u)_{jl}u_l)^{1/2}+C\sum_AM(u)^{AA}\nonumber,
\end{eqnarray}
then diagonalizing  the matrix $(M(u)^{ij}M(u)_{ik}M(u)_{jl})$ at
point $P$ and using the positivity of $M(u)$, we have
\begin{eqnarray}
M(u)^{ij}M(u)_{ik}u_kM(u)_{jl}u_l\ \ \leqq \ \
C\sum_iM(u)^{ii}(\sum_i M(u)_{ii})^2.
\end{eqnarray}
By the above two inequalities,
\begin{eqnarray}
&&M(u)^{Ai}\eta_Au_{ij}u_j\\
&\leqq&
C(M(u)^{AB}\eta_A\eta_B)^{1/2}(\sum_iM(u)^{ii})^{1/2}\sum_iM(u)_{ii}+C\sum_AM(u)^{AA}\nonumber\\
&\leqq&
\frac{\varepsilon_4}{2}\sum_AM(u)^{AA}\sum_iM(u)_{ii}+CM(u)^{AB}\eta_A\eta_B\sum_AM(u)_{AA}+C\sum_AM(u)^{AA}\nonumber.
\end{eqnarray}
Similarly, we have an inequality for the term
$-M(u)^{A\alpha}\eta_Au_{\alpha\beta}u_{\beta}$. Now combining
(4.14), (4.15) and (4.18),
\begin{eqnarray}
\ \ \ \ \ \ \ \ \sum_CM(u)^{AB}\eta_A u_{CCB} &\leqq&
\sum_CM(u)^{AB}\eta_AM(u)_{BCC}+\varepsilon_4\sum_AM(u)^{AA}\sum_AM(u)_{AA}\\
&&+CM(u)^{AB}\eta_A\eta_B\sum_A M(u)_{AA}+C\sum_AM(u)^{AA}\nonumber.
\end{eqnarray}
It will be used in the later. By the Ricci identity, we get
\begin{eqnarray}
  u_{CCAB}&=& u_{ABCC}-u_{EB}R^E_{ACC}-2u_{EC}R^E_{BCA}-u_{AE}R^E_{BCC}.
\end{eqnarray}
Then by (2.9),
\begin{eqnarray}
&&\sum_CM(u)^{AB}u_{CCAB}\\
&=&\sum_CM(u)^{AB}u_{ABCC}-2\sum_CM(u)^{AB}u_{EC}R^E_{BCA}-2\sum_CM(u)^{AB}u_{AE}R^E_{BCC}\nonumber\\
&\geqq&\sum_CM(u)^{AB}u_{ABCC}-C\sum_AM(u)^{AA}\sum_AM(u)_{AA}-C\sum_AM(u)^{AA}-C\nonumber\\
&=&
\sum_CM(u)^{AB}M(u)_{ABCC}+2\sum_CM(u)^{ij}u_{iC}u_{jC}+2\sum_CM(u)^{ij}u_iu_{jCC}\nonumber\\
&&-2\sum_CM(u)^{\alpha\beta}u_{\alpha C}u_{\beta C}-2\sum_CM(u)^{\alpha\beta}u_{\alpha}u_{\beta
CC}
-C\sum_AM(u)^{AA}\sum_AM(u)_{AA}\nonumber\\
&&-C\sum_AM(u)^{AA}-C\nonumber,
\end{eqnarray}
and
\begin{eqnarray}
&&\sum_CM(u)^{ij}u_{iC}u_{jC}-\sum_CM(u)^{\alpha\beta}u_{\alpha C}u_{\beta C}\\
&=&\sum_CM(u)^{iB}u_{iC}u_{BC}-\sum_CM(u)^{A\beta}u_{AC}u_{\beta C}\nonumber\\
&=&\sum_CM(u)^{AB}u_{AC}u_{BC}-2\sum_CM(u)^{A\beta}u_{AC}u_{\beta
C}\nonumber\\
&=& \sum_CM(u)^{AB}u_{AC}M(u)_{BC}
+\sum_jM(u)^{iA}(u_iu_j+\delta_{ij})u_{Aj}\nonumber\\
&&-\sum_{\beta}M(u)^{\alpha A}(u_{\alpha}u_{\beta}+\delta_{\alpha\beta})u_{A\beta}
-2\sum_CM(u)^{A\beta}M(u)_{AC}u_{\beta C}\nonumber\\
&&-2\sum_jM(u)^{i\beta}(u_iu_j+\delta_{ij})u_{\beta j}+2\sum_{\gamma}M(u)^{\alpha\beta}(u_{\alpha}u_{\gamma}+\delta_{\alpha\gamma})u_{\beta\gamma}\nonumber\\
&=& \Delta u+2\sum_j
M(u)^{ik}(u_iu_j+\delta_{ij})u_{kj}+2\sum_{\gamma}M(u)^{\alpha\beta}(u_{\alpha}u_{\gamma}+\delta_{\alpha\gamma})u_{\beta\gamma}\nonumber\\
&&-\sum_jM(u)^{iA}(u_iu_j+\delta_{ij})u_{Aj}-\sum_{\beta}M(u)^{\alpha A}(u_{\alpha}u_{\beta}+\delta_{\alpha\beta})u_{A\beta}
-2\sum_{\alpha}u_{\alpha\alpha}\nonumber\\
&\geqq&
-\frac{\varepsilon_5}{2}\sum_{A,B}u_{AB}^2+2\sum_jM(u)^{ik}M(u)_{kj}(u_iu_j+\delta_{ij})\nonumber\\
&&+2\sum_{\gamma}M(u)^{\alpha\beta}M(u)_{\beta\gamma}(u_{\alpha}u_{\gamma}+\delta_{\alpha\gamma})-C\sum_AM(u)^{AA}-C\nonumber,
\end{eqnarray}
where
$\varepsilon_5=\dfrac{(m+n+2)e^{-2C_3}}{6(1+e^{2C_3})^2(1+C_3^2)^2}$
and we have used the inequality
\begin{eqnarray}
&&-\sum_jM(u)^{iA}(u_iu_j+\delta_{ij})u_{Aj}\\
&=&-\sum_jM(u)^{iA}(u_iu_j+\delta_{ij})M(u)_{Aj}-\sum_jM(u)^{ik}(u_iu_j+\delta_{ij})
(u_ku_j+\delta_{jk})\nonumber\\
&\geqq&-C\sum_AM(u)^{AA}-C\nonumber,
\end{eqnarray}
and a similar inequality for term $-\sum_{\beta}M(u)^{\alpha
A}(u_{\alpha}u_{\beta}+\delta_{\alpha\beta})u_{A\beta}$. Now by
(4.22),
\begin{eqnarray}
&&2(\sum_CM(u)^{ij}u_{iC}u_{jC}-\sum_CM(u)^{\alpha\beta}u_{\alpha
C}u_{\beta C})\\ &\geqq& -\varepsilon_5\sum_{A,B}u_{AB}^2-C\sum_AM(u)^{AA}\sum_AM(u)_{AA}-C\sum_AM(u)^{AA}-C\nonumber.
\end{eqnarray}
Now the term $\sum_{A}M(u)^{AA}\sum_AM(u)_{AA}$ appears, which is
one of our main difficulties. By the Ricci identity, (4.21) becomes
\begin{eqnarray}
&&\sum_CM(u)^{AB}u_{CCAB}\\
&\geqq& \sum_CM(u)^{AB}M(u)_{ABCC}
+2(\sum_CM(u)^{ij}u_iu_{CCj}\nonumber-\sum_CM(u)^{\alpha\beta}u_{\alpha}u_{CC\beta})\\
&&-\varepsilon_5\sum_{A,B}u_{AB}^2-C\sum_AM(u)^{AA}\sum_AM(u)_{AA}-C\sum_AM(u)^{AA}-C\nonumber.
\end{eqnarray}
At point $P$, by (4.3),(4.10),(4.12),(4.13),(4.25) we have
\begin{eqnarray}
&&L(\phi)\\ &\geqq& e^{-C_6\eta}[C_6^2M(u)^{AB}\eta_A\eta_B\Delta u-C_6M(u)^{AB}\eta_{AB}\Delta u
-2C_6M(u)^{AB}\eta_A\sum_Cu_{CCB}\nonumber\\
&& +\sum_CM(u)^{AB}M(u)_{ABCC} +2C_6M(u)^{ij}u_i\eta_{j}\Delta u-2C_6M(u)^{\alpha\beta}u_{\alpha}\eta_{\beta}\Delta
u\nonumber\\
&&-\varepsilon_5\sum_{A,B}u_{AB}^2-C\sum_AM(u)^{AA}\sum_AM(u)_{AA}-C\sum_AM(u)^{AA}-C]\nonumber\\
&&+M(u)^{AB}\zeta_{AB}-2M(u)^{ij}u_i\zeta_j+2M(u)^{\alpha\beta}u_{\alpha}\zeta_{\beta}\nonumber\\
&=& e^{-C_6\eta}\{C_6^2M(u)^{AB}\eta_A\eta_B\Delta
u-2C_6M(u)^{AB}\eta_A\sum_Cu_{CCB}+\sum_CM(u)^{AB}M(u)_{ABCC}\nonumber\\
&&
+C_6[M(u)^{AB}(-\eta)_{AB}+2M(u)^{ij}u_i\eta_{j}-2M(u)^{\alpha\beta}u_{\alpha}\eta_{\beta}]\sum_AM(u)_{AA}\nonumber\\
&&+C_6[M(u)^{AB}(-\eta)_{AB}+2M(u)^{ij}u_i\eta_{j}-2M(u)^{\alpha\beta}u_{\alpha}\eta_{\beta}](|\nabla_x u|^2-|\nabla_y u|^2+m-n)\nonumber\\
&&-\varepsilon_5\sum_{A,B}u_{AB}^2-C\sum_AM(u)^{AA}\sum_AM(u)_{AA}-C\sum_AM(u)^{AA}-C\}\nonumber\\
&&+M(u)^{AB}\zeta_{AB}-2M(u)^{ij}u_i\zeta_j+2M(u)^{\alpha\beta}u_{\alpha}\zeta_{\beta}\nonumber
\end{eqnarray}
Now using the first equality of (4.4), the bounds on $\psi$, on
$x_1$, and on the $C^1$ norm of $u$, we have
\begin{eqnarray}
M(u)^{AB}(-\eta)_{AB}+2M(u)^{ij}u_i\eta_{j}-2M(u)^{\alpha\beta}u_{\alpha}\eta_{\beta}&\leqq&C\sum_{A}M(u)^{AA}+C.
\end{eqnarray}
Combing the above two inequalities and (4.10), we have
\begin{eqnarray}
&&L(\phi)\\
&\geqq& e^{-C_6\eta}\{C_6^2M(u)^{AB}\eta_A\eta_B\Delta u-2C_6M(u)^{AB}\eta_A\sum_Cu_{CCB}+\sum_CM(u)^{AB}M(u)_{ABCC}\nonumber\\
&&+(4\varepsilon_4C_6-C)\sum_AM(u)^{AA}\sum_AM(u)_{AA}-\varepsilon_5\sum_{A,B}u_{AB}^2-(CC_6+C)\sum_AM(u)^{AA}\nonumber\\
&&-C(C_6)\}+4\varepsilon_4\sum_AM(u)^{AA}-C\nonumber.
\end{eqnarray}
Here $C(C_6)$ is a positive constant depending on $\psi,m,n$, the
$C^1$ norm of $u$ and also $C_6$. Now we use equation. Take
logarithm of (1.2) and differentiate it:
\begin{eqnarray}
M(u)^{AB}M(u)_{ABC}&=& [\ln(K)]_C+\frac{m+n+2}{2}[\ln(f)]_C,\\
M(u)^{AB}M(u)_{ABCC}
&=&M(u)^{AA^{'}}M(u)^{BB^{'}}M(u)_{ABC}M(u)_{A^{'}B^{'}C}\\ &&+[\ln(
K)]_{CC} +\frac{m+n+2}{2}[\ln(f)]_{CC}\nonumber.
\end{eqnarray}
We choose a local orthonormal frame at $P$ such that the matrix
$M(u)$ is diagonal at $P$. Then by (4.30),
\begin{eqnarray}
\ \ \ \ \sum_CM(u)^{AB}M(u)_{ABCC}
&\geqq&\sum_{A,B,C}\frac{M(u)^2_{ABC}}{M(u)_{AA}M(u)_{BB}}-\frac{m+n+2}{2}\sum_C[\ln(f)]^2_C\\
&&
+\frac{m+n+2}{2}\sum_C\frac{(f)_{CC}}{f}-C\nonumber,
\end{eqnarray}
where we used the bound of $K$. Obviously,
\begin{eqnarray}
&&\sum_{A,B,C}\frac{M(u)^2_{ABC}}{M(u)_{AA}M(u)_{BB}}\\
&\geqq&\sum_{B\neq C}\frac{M(u)^2_{CBC}}{M(u)_{CC}M(u)_{BB}}
+\sum_{A\neq C}\frac{M(u)^2_{ACC}}{M(u)_{AA}M(u)_{CC}}+\sum_{A,C}\frac{M(u)^2_{AAC}}{M(u)^2_{AA}}\nonumber.
\end{eqnarray}
Then by (4.29) and Schwarz inequality,
\begin{eqnarray}
&&\frac{m+n+2}{4}\sum_C[\ln(f)]^2_C \ \ =\ \ \frac{1}{m+n+2}\sum_C[\sum_{A}\frac{M(u)_{AAC}}{M(u)_{AA}}-(\ln(K))_C]^2\\
&\leqq& \frac{1}{m+n+2}(1+\frac{1}{m+n+1})\sum_C[\sum_A\frac{M(u)_{AAC}}{M(u)_{AA}}]^2+C
\nonumber\\
&\leqq& \frac{m+n}{m+n+1}\sum_{A,C}[\frac{M(u)_{AAC}}{M(u)_{AA}}]^2+C\nonumber.
\end{eqnarray}
Now by (2.1)
\begin{eqnarray}
(f)_C=f_ru_C+2[f_p\sum_{i}u_iu_{iC}+f_q\sum_{\alpha }u_{\alpha}u_{\alpha C}],
\end{eqnarray}
then  by the Schwarz inequality and $f_p$, $f_q$ $>$ $0$,
\begin{eqnarray}
&&\frac{m+n+2}{4}\sum_C[\ln(f)]^2_C\ \ = \ \
\frac{m+n+2}{4}\sum_C[\frac{(f)_C}{f}]^2\\
&\leqq&
C+\varepsilon_5\sum_{AB}u_{AB}^2+\frac{m+n+2}{f^2}\sum_C[f_p\sum_{i}u_iu_{iC}+f_q\sum_{\alpha}u_{\alpha}u_{\alpha
C}]^2\nonumber\\
&\leqq&C+\varepsilon_5\sum_{AB}u_{AB}^2+\frac{m+n+2}{f^2}[f_p\sum_iu^2_i+f_q\sum_{\alpha}u^2_{\alpha}][f_p\sum_{i,C}u^2_{iC}
+f_q\sum_{\alpha, C}u^2_{\alpha C}]\nonumber.
\end{eqnarray}
Now by (2.1) and Proposition 3.3,
\begin{eqnarray}
&&\frac{m+n+2}{f^2}[f-f_p\sum_iu^2_i-f_q\sum_{\alpha}u^2_{\alpha}][f_p\sum_{i,C}u^2_{iC}+f_q\sum_{\alpha, C}u^2_{\alpha
C}]\\
&=&\frac{m+n+2}{(1+e^{2u}+ |\nabla_x u|^2+e^{2u}|\nabla_y u|^2)^2}[(1+e^{2u})\sum_{i,C}u^2_{iC}+e^{2u}(1+e^{2u})\sum_{\alpha, C}u^2_{\alpha
C}]\nonumber\\
&\geqq&\frac{(m+n+2)e^{-2C_3}}{(1+e^{2C_3})^2(1+C_3^2)^2}\sum_{A,B}u^2_{AB}\
\ \geqq \ \ 6\varepsilon_5 \sum_{A,B}u^2_{AB}\nonumber.
\end{eqnarray}
Now (4.35) becomes
\begin{eqnarray}
&&\frac{m+n+2}{4}\sum_C[\ln(f)]^2_C\\
&\leqq&C-5\varepsilon_5\sum_{AB}u_{AB}^2+\frac{m+n+2}{f}[f_p\sum_{i,C}u^2_{iC}+f_q\sum_{\alpha, C}u^2_{\alpha C}]\nonumber.
\end{eqnarray}
By (4.34),
\begin{eqnarray}
&&\frac{m+n+2}{2}\sum_C\frac{(f)_{CC}}{f}\\
&=&\frac{m+n+2}{2f}\sum_C(f_ru_C)_C+\frac{m+n+2}{f}\sum_{i,C}(f_p)_Cu_iu_{iC}+\frac{m+n+2}{f}\sum_{\alpha,C}(f_q)_Cu_{\alpha}u_{\alpha
C}\nonumber\\
&&+\frac{m+n+2}{f}[f_p\sum_{i,C}u^2_{iC}+f_q\sum_{\alpha,
C}u^2_{\alpha
C}]+\frac{m+n+2}{f}[f_p\sum_{i,C}u_iu_{iCC}+f_q\sum_{\alpha,
C}u_{\alpha}u_{\alpha
CC}]\nonumber\\
&\geqq&-\varepsilon_5\sum_{AB}u^2_{AB}-C+\frac{m+n+2}{f}[f_p\sum_{i,C}u^2_{iC}+f_q\sum_{\alpha, C}u^2_{\alpha
C}]\nonumber\\
&&+\frac{m+n+2}{f}[f_p\sum_{i,C}u_iu_{CCi}+f_q\sum_{\alpha, C}u_{\alpha}u_{CC \alpha}]\nonumber.
\end{eqnarray}
By (4.12),
\begin{eqnarray}
\frac{m+n+2}{f}[f_p\sum_{i,C}u_iu_{CCi}+f_q\sum_{\alpha,
C}u_{\alpha}u_{CC \alpha}]&\geqq&-C(C_6)-\varepsilon_5\sum_{A,B}u_{AB}^2,
\end{eqnarray}
where $C(C_6)$ is a positive constant depending on an undetermined
constant $C_6$. Then by (4.33), (4.37), (4.38), (4.39),
\begin{eqnarray}
&&-\frac{m+n+2}{2}\sum_C[\ln(f)]^2_C+\frac{m+n+2}{2}\sum_C\frac{(f)_{CC}}{f}\\
&\geqq&-\frac{m+n}{m+n+1}\sum_{A,C}[\frac{M(u)_{AAC}}{M(u)_{AA}}]^2+3\varepsilon_5\sum_{A,B}u_{AB}^2-C(C_6)\nonumber.
\end{eqnarray}
Now by (4.31), (4.32), (4.40) and (4.1), $\eta<3$,
\begin{eqnarray}
&&\sum_CM(u)^{AB}M(u)_{ABCC}\\
&\geqq& 2\sum_{A\neq
C}\frac{M(u)^2_{ACC}}{M(u)_{AA}M(u)_{CC}}+\frac{1}{m+n+1}\sum_{A,C}\frac{M(u)^2_{AAC}}{M(u)^2_{AA}}+3\varepsilon_5\sum_{A,B}u_{AB}^2-C(C_6)\nonumber.
\end{eqnarray}
Then by (4.19),(4.28),(4.41),  and Ricci identity, we have
\begin{eqnarray}
&&L(\phi)\\
&\geqq&
e^{-C_6\eta}\{C_6(C_6-\hat{C}_6)M(u)^{AB}\eta_A\eta_B\sum_AM(u)_{AA}
-2C_6\sum_CM(u)^{AB}\eta_AM(u)_{BCC}\nonumber\\&&+(2C_6\varepsilon_4-C)\sum_AM(u)^{AA}\sum_AM(u)_{AA}
+2\sum_{A\neq C}\frac{M(u)^2_{ACC}}{M(u)_{AA}M(u)_{CC}}\nonumber\\&&
+\frac{1}{m+n+1}\sum_{A,C}\frac{M(u)^2_{AAC}}{M(u)^2_{AA}}+2\varepsilon_5\sum_{A,B}u_{AB}^2
-(CC_6^2+CC_6+C)\sum_AM(u)^{AA}\nonumber\\&&-C(C_6)\}
+4\varepsilon_4\sum_AM(u)^{AA}-C\nonumber.
\end{eqnarray}
Now take $C_6>\hat{C}_6$. Since
\begin{eqnarray}
&&2C_6\sum_{A\neq C}\frac{\eta_AM(u)_{ACC}}{M(u)_{AA}}\ \ = \ \ C_6\sum_{A\neq
C}\frac{2(M(u)_{CC}\eta_A)M(u)_{ACC}}{M(u)_{AA}M(u)_{CC}}\\
&\leqq& C_6 \sum_{A\neq
C}\frac{(C_6-\hat{C}_6)M(u)^2_{CC}\eta^2_A+\dfrac{M(u)^2_{ACC}}{C_6-\hat{C}_6}}{M(u)_{AA}M(u)_{CC}}\nonumber\\
&\leqq& C_6(C_6-\hat{C}_6)\sum_{A,
C}\frac{\eta^2_AM(u)_{CC}}{M(u)_{AA}}+\frac{C_6}{C_6-\hat{C}_6}\sum_{A\neq
C}\frac{M(u)^2_{ACC}}{M(u)_{AA}M(u)_{CC}}\nonumber,
\end{eqnarray}
and
\begin{eqnarray}
\ \ \ \ \  2C_6\sum_A \frac{\eta_AM(u)_{AAA}}{M(u)_{AA}} &\leqq&
\frac{1}{(m+n+1)}\sum_A \frac{M(u)^2_{AAA}}{M(u)^2_{AA}}
+(m+n+1)C^2_6\sum_A\eta_A^2,
\end{eqnarray}
then
\begin{eqnarray}
&&C_6(C_6-\hat{C}_6)M(u)^{AB}\eta_A\eta_B\sum_AM(u)_{AA}-2C_6\sum_CM(u)^{AB}\eta_AM(u)_{BCC}\\
&&+2\sum_{A\neq
C}\frac{M(u)^2_{ACC}}{M(u)_{AA}M(u)_{CC}}+\frac{1}{m+n+1}\sum_{A,C}\frac{M(u)^2_{AAC}}{M(u)^2_{AA}}\nonumber\\
&\geqq& \big( 2-\frac{C_6}{C_6-\hat{C}_6}\big)\sum_{A\neq
C}\frac{M(u)^2_{ACC}}{M(u)_{AA}M(u)_{CC}}-CC_6^2\nonumber.
\end{eqnarray}
Now by (4.1), (4.45) , $x_1\leqq 1$, we have
\begin{eqnarray}
&&L(\phi)\\
&\geqq& e^{-C_6\eta}\{\big(
2-\frac{C_6}{C_6-\hat{C}_6}\big)\sum_{A\neq
C}\frac{M(u)^2_{ACC}}{M(u)_{AA}M(u)_{CC}}-(\hat{C}_8C_6^2+\hat{C}_9)\sum_AM(u)^{AA}\nonumber\\
&&+(2C_6\varepsilon_4-\hat{C}_7)\sum_AM(u)^{AA}\sum_AM(u)_{AA}+\varepsilon_5\sum_{A,B}u_{AB}^2\}+4\varepsilon_4\sum_AM(u)^{AA}
-C(C_6)\nonumber.
\end{eqnarray}
We take $C_6$ big enough to satisfy
\begin{eqnarray}
&&2-\frac{C_6}{C_6-\hat{C}_6}\ \ >\ \ 0,\ \ \ \
2C_6\varepsilon_4-\hat{C}_7\ \ > \ \ 0, \ \ \text{and} \ \
4\varepsilon_4-e^{-C_6}(\hat{C}_8C_6^2+\hat{C}_9)\ \
>\ \ 0.
\end{eqnarray}
Now by (4.47) and $\eta \geqq 1$, (4.46) implies that
$\sum_{A,B}u^2_{AB}$ is bounded  at point $P$. So the function
$\phi$ has a uniform upper bound, and by the positivity of matrix
$M(u)$, there is a constant $C_7$ depending on $K$, $\psi$, $m$,
$n$, $\partial\Omega$ such that
\begin{eqnarray}
\sum_{A,B}|u_{AB}|^2&\leqq&C_7.
\end{eqnarray}
This gives the interior $C^2$ estimate. Here we generalize the idea
of [Y] to deal with the $C^3$ term. In order to obtain the $C^2$
estimate, now we only need the estimate on the boundary.
\section{$C^2$ estimate on the boundary }\noindent
Let $P$ be on the boundary $\partial \Omega_x\times S^n$,
$P=(P_x,P_y)$, and $\Omega^{\delta}(P)=\Omega^{\delta}_x(P_x)\times
B^{\delta}_y(P_y)$. Here  $\Omega_x^{\delta}(P_x)=\Omega_x\cap
B^{\delta}_x(P_x)$, and $B^{\delta}_x(P_x),B^{\delta}_y(P_y)$ are
$\delta$ geodesic sphere neighborhoods of $S^m$ and $ S^n$ centered
at $P_x$ and $P_y$ respectively. Since $\Omega$ is a PHC-domain, for
sufficiently small $\delta$, we can find a frame
$\{e_1,\cdots,e_{m+n}\}$ on $\bar{\Omega}^{\delta}(P)$ such that:
$e_m$ is the outer normal direction on $\partial\Omega^{\delta}(P)$;
the previous $m-1$ ones are tangent vectors of $\partial\Omega_x$;
the last $n$ ones are tangent vectors of $S^n$. By compactness of
$\partial\Omega_x$, we can take $\delta$ independent from boundary
points. (The proof is similar to Lebesgue's Covering Lemma.) Taking
$\delta <1$ sufficiently small, we consider a local function on
$\bar{\Omega}^{\delta}(P)$,
\begin{eqnarray}
\phi&=&(u-\psi)_C+v,
\end{eqnarray}
where $\psi\in$ ABF($\tau,K$), and $v$ is an undetermined function
which we will give explicitly in the following. If $\phi$ attains
its maximum value at some point $Q$ in $\Omega^{\delta}(P)$. Then at
point $Q$,
\begin{eqnarray}
(u-\psi)_{CA}+v_A&=&0.
\end{eqnarray}
By (3.23) and Ricci identity,
\begin{eqnarray}
L\phi & = & M(u)^{AB}(u-\psi)_{CAB}+M(u)^{AB}v_{AB}\\
&\geqq&M(u)^{AB}u_{ABC}+M(u)^{AB}v_{AB}-C(1+\sum_AM(u)^{AA})\nonumber.
\end{eqnarray}
For any function $\xi$ define an operator
\begin{eqnarray}
\tilde{L}(\xi)&=&M(u)^{AB}\xi_{AB}-2M(u)^{ij}u_i\xi_{j}+2M(u)^{\alpha\beta}u_{\alpha}\xi_{\beta }\\
&&-\frac{m+n+2}{f}[f_p\sum_iu_i\xi_i+f_q\sum_{\alpha}u_{\alpha}\xi_{\alpha}]\nonumber.
\end{eqnarray}
Then by (4.29),(5.2) and (5.4), (5.3) becomes
\begin{eqnarray}
&&L\phi\\
&\geqq& 2M(u)^{ij}u_iu_{jC}-2M(u)^{\alpha\beta}u_{\alpha}u_{\beta C}+\frac{m+n+2}{f}(f_p\sum_iu_iu_{jC}+f_q\sum_{\alpha}u_{\alpha}u_{\alpha
C})\nonumber\\
&&+M(u)^{AB}v_{AB}-C(1+\sum_AM(u)^{AA})\nonumber\\
&\geqq& \tilde{L}v-C(1+\sum_{A}M(u)^{AA})\nonumber.
\end{eqnarray}
Now we define a vector field in $\mathbb{R}^{m+1}$. For
$0<\epsilon<1/4$ and $P_x\in \partial\Omega_x$, let
\begin{eqnarray}
\vec{\chi}(P_x)&=& -e_m(P_x)+\epsilon\vec{\gamma}(P_x).
\end{eqnarray}
Here we use $\vec{\gamma}(\cdot)$ and $e_m(\cdot)$ to denote the
vectors at point "$\cdot$" in $\mathbb{R}^{m+1}$, and $\vec{\gamma}$
is defined in (1.1). We choose a coordinate system for
$\mathbb{R}^{m+1}$ with first coordinate axis given by
$\vec{\chi}(P_x)$. Denote by $<\cdot,\cdot>_m$ and $|\cdot|_m$ the
inner product and corresponding norm of $\mathbb{R}^{m+1}$. Let
$H=\{\vec{\gamma}\in S^m;<-e_m(P_x),\vec{\gamma}>_m=0\}$ be a
totally geodesic submanifold of $S^m$. Since $\Omega_x $ is a
strictly infinitesimally convex domain, it is a strictly locally
convex domain (see [S]). Since the exponential map of $S^m$ takes
the subspace $T_{P_x}(\partial\Omega_x) $ onto $H$, we have $H\cap
\bar{\Omega}_x=\{P_x\}$. (If there is another point $\tilde{P}_x\in
H\cap \bar{\Omega}_x$, the minimal geodesic curve connected $P_x$
and $\tilde{P}_x$ is contained in $H\cap\bar{\Omega}_x$ which
contradicts the strictly locally convexity of $\partial\Omega_x$.)
This means that for any $Q'_x\in \Omega^{\delta}_x$,
\begin{eqnarray}
<\vec{\gamma}(Q_x')-\vec{\gamma}(P_x),-e_m(P_x)>_m&>&0.
\end{eqnarray}
Now define $S_{\delta'}=\{Q'_x\in
S^m;|\vec{\gamma}(Q'_x)-\vec{\gamma}(P_x)|_m=\delta'\}$. For
$Q'_x\in \Omega^{\delta}_x(P_x)\cap S_{\delta'}$ where
$\delta'<2\sin(\delta/2)$ (guaranteeing that $\Omega^{\delta}_x\cap
S_{\delta'}$ is non-empty), we have
\begin{eqnarray}
<\vec{\gamma}(Q'_x)-\vec{\gamma}(P_x),\vec{\gamma}(P_x)>_m &=&
-(\delta')^2/2,
\end{eqnarray}
where we used
$|\vec{\gamma}(Q'_x)-\vec{\gamma}(P_x)|_m^2=(\delta')^2 $ and
$|\vec{\gamma}(Q'_x)|_m=|\vec{\gamma}(P_x)|_m=1 $. Now denote
$a_i(Q'_x)=<\vec{\gamma}(Q'_x)-\vec{\gamma}(P_x),e_i(P_x)>_m$. Then
by (5.8), $\sum_i^ma_i^2(Q'_x)+(\delta')^4/4=(\delta')^2$. By (5.7),
we have $-a_m(Q'_x)>0$. Now further assume $Q'_x$ is the minimum
value point of function $-a_m$. Then
$a_i^2(Q'_x)<(\delta')^2-(\delta')^4/4$ for $i\neq m$. We can take a
vector $\vec{b}$ in $S_{\delta'}$ defined by
$\vec{b}=\sum_ib_ie_i(P_x)+[1-(\delta')^2/2]\vec{\gamma}(P_x)$ such
that for $i\neq m$, $|b_i|>|a_i(Q'_x)|$; and $-b_m<-a_m(Q'_x)$. By
the openness of  $\Omega^{\delta}_x(P_x)$, the point corresponding
to $\vec{b}$ is in $\Omega^{\delta}_x(P_x)$ if we further require
that $|b_i-a_i(Q'_x)|$ is sufficiently small. This is a
contradiction. Hence the minimal value of function $-a_m$ in
$\bar{\Omega}^{\delta}_x(P_x)\cap S_{\delta'}$ only occurs in
$\partial\Omega^{\delta}_x(P_x)\cap S_{\delta'}$. So for any
$Q'_x\in \Omega^{\delta}_x(P_x)\cap S_{\delta'}$, by  (5.8) hold for
any point in $S_{\delta'}$ and (5.6),
\begin{eqnarray}
\ \ <\vec{\gamma}(Q'_x)-\vec{\gamma}(P_x),\vec{\chi}(P_x)>_m &\geqq&
\inf_{Q''_x\in
\partial\Omega^{\delta}_x(P_x)\cap S_{\delta'}
}<\vec{\gamma}(Q''_x)-\vec{\gamma}(P_x),\vec{\chi}(P_x)>_m.
\end{eqnarray}
Let $\nabla,\bar{\nabla}$ be the Levi-Civita connections of $S^m$
and $\mathbb{R}^{m+1}$. Since $\partial\Omega_x$ is a strictly
infinitesimally convex hypersurface in $S^m$, its second fundamental
tensor is positive definite everywhere.  Then  for $i,j\neq m$, on
$\partial\Omega^{\delta}_x(P_x)$ the order $m-1$ matrix
\begin{eqnarray}
<\bar{\nabla}_{e_i}\bar{\nabla}_{e_j}\vec{\gamma},-e_m>_m &=&
<\nabla_{e_i}e_j,-e_m>_m
\end{eqnarray}
has a positive uniform (independent from the choice of boundary
points) lower bound, by the compactness of $\partial\Omega_x$. By
(5.6), we can take a uniformly sufficiently small $\epsilon$ such
that the matrix
$(<\bar{\nabla}_{e_i}\bar{\nabla}_{e_j}\vec{\gamma},\vec{\chi}>_m)_{(m-1)\times
(m-1)}(P_x)$ also has a uniform  positive lower bound. Now by the
Taylor expansion of $\partial\Omega_x$ near $P_x$ with the frame
$e_1(P_x),\cdots,e_{m-1}(P_x)$, we find that the right hand side of
(5.9) is non-negative for sufficiently small $\delta$. Moreover the
choice of $\delta$
 is independent from the boundary point, as follow from: (i)
there is a uniform $\delta$ such that:
$\partial\Omega^{\delta}_x(P_x)$ can be parameterized, and the
tangent vectors along the parameterized curves at $P_x$ are
$e_1(P_x),\cdots,e_{m-1}(P_x)$; (ii) for the function
$<\vec{\gamma}(Q_x'')-\vec{\gamma}(P_x),\vec{\chi}(P_x)>_m$ where
 $Q''_x\in \partial\Omega^{\delta}_x(P_x)$, the second order term at $P_x$ has a uniform lower
 bound,
 and higher than second order terms at $P_x$ have uniform upper bound. By the arbitrary choice
of $\delta'$, we have
 $<\vec{\gamma}(Q_x')-\vec{\gamma}(P_x),\vec{\chi}(P_x)>_m\geqq
0$ for $Q_x'\in \Omega^{\delta}_x(P_x)$£¬which implies in
$\Omega^{\delta}_x(P_x)$ that
\begin{eqnarray}
x_1& \geqq & x_1(P_x)\ \ = \ \ \frac{\epsilon}{(1+\epsilon^2)^{1/2}}
\ \ > \ \ 0 .
\end{eqnarray}
Denote $\theta=x_1(P_x)$, which is  a constant only depending on
$\partial\Omega_x$.  Now obviously for $i\neq m$,
\begin{eqnarray}
(x_1)_i(P_x)&=&0,\ \ \text{and} \ \  |(x_1)_m(P_x)|\ \ \geqq\ \
\frac{3}{4}.
\end{eqnarray}
Now we let
\begin{eqnarray}
w&=&e^{-C_8(u-\psi)}+e^{-C_{9}C_8(x_1-\theta)}.
\end{eqnarray}
Here $C_8$, $C_9$ are two constants which we will determine in the
following. Using the similar tricks of  (4.4) to (4.9), we have
\begin{eqnarray}
&&\ \ \tilde{L}(e^{-C_8(u-\psi)})\\
&\geqq&e^{-C_9(u-\psi)}[-\hat{C}_{10}C_8+(3\varepsilon_3C_8-\hat{C}_{11})\sum_AM(u)^{AA}-\hat{C}_{12}C_8\sum_iM(u)^{ii}
]\nonumber,
\end{eqnarray}
and in $\Omega^{\delta}_x(P_x)$,
\begin{eqnarray}
&&\tilde{L}(e^{-C_{9}C_8(x_1-\theta)})\\
&=&e^{-C_{9}C_8(x_1-\theta)}\{(C_{9}C_8)^2M(u)^{ij}(x_1)_i(x_1)_j-C_{9}C_8M(u)^{ij}(x_1)_{ij}\nonumber\\
&&+2C_{9}C_{8}M(u)^{ij}u_i(x_1)_j+\frac{(m+n+2)C_{9}C_8}{f}f_p\sum_iu_i(x_1)_i\}\nonumber\\
&\geqq&e^{-C_{9}C_8(x_1-\theta)}\{\frac{(C_{9}C_8)^2}{2}M(u)^{ij}(x_1)_i(x_1)_j
+(C_{9}C_8\theta-C)\sum_iM(u)^{ii}-CC_{9}C_8\}\nonumber,
\end{eqnarray}
where we used $(x_1)_{ij}+x_1\delta_{ij}=0$, $x_1\geqq\theta$ and
the inequality
\begin{eqnarray}
|2C_{9}C_8M(u)^{ij}u_i(x_1)_j|&\leqq&\frac{(C_{9}C_8)^2}{2}M(u)^{ij}(x_1)_i(x_1)_j+C\sum_iM(u)^{ii}.
\end{eqnarray}
Then for any point $Q'$ $\in$ $\Omega^{\delta}(P)$,
$Q'=(Q_x',Q_y')$, we have
\begin{eqnarray}
|(x_1)_i(P_x)-(x_1)_i(Q_x')|&\leqq&|\nabla(x_1)_i|dist_x(P_x,Q_x')\
\ \leqq  \ \  \hat{C}_{13}\delta,
\end{eqnarray}
where $\hat{C}_{13}$ is an absolute positive constant and
$dist_x(\cdot,\cdot)$ is the distance function of $S^m$. Assume
$C_{9},C_8>1$. Then we choose $\delta$ such that
\begin{eqnarray}
\delta&\leqq& \frac{1}{6\hat{C}_{13}(C_{9}C_8)^2}.
\end{eqnarray}
By  (5.12), (5.17), we have for $i\neq m$
\begin{eqnarray}
|(x_1)_i(Q_x')|\ \ \leqq\ \ \frac{1}{6(C_{9}C_8)^2}\ \ <\ \ 1,\ \
\text{and} \ \ \frac{1}{2} \ \ < \ \ |(x_1)_m(Q^{'}_x)|.
\end{eqnarray}
Then in $\Omega^{\delta}(P)$, by the positivity of $M(u)$ and
(5.16),
\begin{eqnarray}
&&M(u)^{ij}(x_1)_i(x_1)_j\\
&=& M(u)^{mm}(x_1)_m^2+\sum_{i,j\neq
m}M(u)^{ij}(x_1)_i(x_1)_j+ 2M(u)^{mi}(x_1)_i(x_1)_m\nonumber\\
& \geqq&
\frac{1}{4}M(u)^{mm}-\frac{1}{6(C_{9}C_8)^2}\sum_iM(u)^{ii}-\sum_i(M(u)^{mm}+M(u)^{ii})|(x_1)_i||(x_1)_m|\nonumber\\
&\geqq&\frac{1}{4}M(u)^{mm}-\frac{C}{(C_{9}C_8)^2}\sum_iM(u)^{ii}\nonumber.
\end{eqnarray}
Now (5.12) becomes
\begin{eqnarray}
&&\ \ \ \ \tilde{L}(e^{-C_{9}C_8(x_1-\theta)})\\
&\geqq&e^{-C_{9}C_8(x_1-\theta)}[\frac{(C_{9}C_8)^2}{8}M(u)^{mm}+(C_{9}C_8\theta-\hat{C}_{14})\sum_iM(u)^{ii}-\hat{C}_{15}C_{9}C_8]\nonumber.
\end{eqnarray}
Since at point $P$, $u-\psi-C_{9}(x_1-\theta)=0$, we know that in
$\Omega^{\delta}(P)$,
\begin{eqnarray}
\ \
 |u-\psi-C_{9}(x_1-\theta)|&\leqq&(|\nabla(u-\psi)|+C_{9}|\nabla(x_1)|)\delta\
\ \leqq\ \ \hat{C}_{16}C_{9}\delta,
\end{eqnarray}
where $\hat{C}_{16}$ is a positive constant depending on $\psi$ and
$C_3$. We further require
\begin{eqnarray}
\delta&\leqq&\frac{1}{C_{9}C_8\hat{C}_{16}}.
\end{eqnarray}
We take $C_{9}=\max\{1,(e\hat{C}_{12}+\hat{C}_{14})/\theta\}$, and
\begin{eqnarray}
C_8&\geqq&\frac{\hat{C}_{11}}{\varepsilon_3}.
\end{eqnarray}
Now by (5.14), (5.21), (5.22) and (5.23),
\begin{eqnarray}
&&\tilde{L}(e^{-C_8(u-\psi)}+e^{-C_{9}C_8(x_1-\theta)})\\
&\geqq&e^{-C_8(u-\psi)}\{-\hat{C}_{10}C_8+(3\varepsilon_3C_8-\hat{C}_{11})\sum_AM(u)^{AA}-\hat{C}_{12}C_8\sum_iM(u)^{ii}\nonumber\\
&&+
e^{C_8[(u-\psi)-C_{9}(x_1-\theta)]}[\frac{(C_{9}C_8)^2}{8}M(u)^{mm}+(C_{9}C_8\theta-\hat{C}_{14})\sum_iM(u)^{ii}-\hat{C}_{15}C_{9}C_8]\}\nonumber\\
&\geqq&e^{-C_8(u-\psi)}[-\hat{C}_{10}C_8+(3\varepsilon_3C_8-\hat{C}_{11})\sum_AM(u)^{AA}-\hat{C}_{12}C_8\sum_iM(u)^{ii}\nonumber\\
&&+
e^{-1}\frac{(C_{9}C_8)^2}{8}M(u)^{mm}+e^{-1}(C_{9}C_8\theta-\hat{C}_{14})\sum_iM(u)^{ii}-e\hat{C}_{15}C_{9}C_8]\nonumber.
\end{eqnarray}
Since $C_8>1$ and by (1.2),(5.13),(5.24) and (5.25),
\begin{eqnarray}
&&\ \ \tilde{L}(w)\\
&\geqq&e^{-C_8(u-\psi)}[2\varepsilon_3C_8\sum_AM(u)^{AA}+e^{-1}\frac{(C_{9}C_8)^2}{8}M(u)^{mm}-(\hat{C}_{10}+e\hat{C}_{15}C_{9})C_8]\nonumber\\
&\geqq&e^{-C_8(u-\psi)}[\varepsilon_3C_8\sum_AM(u)^{AA}+\varepsilon_3C_8\sum_{A\neq
1}\lambda_A+e^{-1}\frac{(C_{9}C_8)^2}{8}\lambda_1-(\hat{C}_{10}+e\hat{C}_{15}C_{9})C_8]\nonumber\\
&\geqq&e^{-C_8(u-\psi)}[\varepsilon_3C_8\sum_AM(u)^{AA}+(C_8)^{\frac{m+n+1}{m+n}}C(\prod_A\lambda_A)^{\frac{1}{m+n}}-(\hat{C}_{10}+e\hat{C}_{15}C_{9})C_8]\nonumber\\
&\geqq&e^{-C_8(u-\psi)}[\varepsilon_3C_8\sum_AM(u)^{AA}+(C_8)^{\frac{m+n+1}{m+n}}\hat{C}_{17}-(\hat{C}_{10}+e\hat{C}_{15}C_{9})C_8]\nonumber,
\end{eqnarray}
where we assume $\lambda_1\leqq\lambda_2 \cdots\leqq\lambda_{m+n}$
to be the positive eigenvalues of the matrix $M(u)^{-1}$, and
$\hat{C}_{17}$ is a positive constant depending on $\varepsilon_3$,
$C_{9}$, $C_3$. Now by (5.24), we take
\begin{eqnarray}
C_{8}&=& \max\{1,
\frac{\hat{C}_{11}}{\varepsilon_3},\frac{(\hat{C}_{10}+e\hat{C}_{15}C_{9})^{m+n}}{\hat{C}_{17}^{m+n}}\}.
\end{eqnarray}
Now we choose $\delta$ sufficiently small and satisfying
(5.18),(5.23), then in $\Omega^{\delta}(P)$ by (5.26),(5.27),
\begin{eqnarray}
\tilde{L}(w)&\geqq& \hat{C}_{18}\sum_AM(u)^{AA}.
\end{eqnarray}
For $Q'\in \Omega^{\delta}(P)$, we let
\begin{eqnarray}
v(Q')&=& C_{10}w(Q')-C_{11}d^2(Q'),
\end{eqnarray}
where $d(\cdot)=dist(P,\cdot)$, $dist(\cdot,\cdot)$ is the distance
function of $S^m\times S^n$, and $C_{10},C_{11}$ are two positive
constants which will be determined in the following. Then by (5.1)
and (5.29), we know that on
$\partial\Omega\cap\bar{\Omega}^{\delta}(P) \setminus{\{P\}}$ (where
$w<2$ and $u|_{\partial\Omega}=\psi$), we have
\begin{eqnarray}
\phi&=& C_{10}w-C_{11}d^2\ \ < \ \  2C_{10}.
\end{eqnarray}
Obviously $\phi(P)=2C_{10}$. Moreover, on $\partial
\Omega^{\delta}\cap \bar{\Omega}$,
\begin{eqnarray}
\phi&\leqq&(u-\psi)_C+2C_{10}-C_{11}d^2\ \ \leqq\ \ \hat{C}_{19}
+2C_{10}-C_{11}\delta^2,
\end{eqnarray}
where $\hat{C}_{19}$ is a positive constant depending on $C_3$,
$\psi$. Now we take $C_{11}=(\hat{C}_{19}+1)/\delta^2$. By (5.30)
and (5.31) on $\partial\Omega^{\delta}(P)\setminus \{P\}$,
\begin{eqnarray}
\phi&<&\phi(P)\ \ =\ \ 2C_{10}.
\end{eqnarray}
We notice that the derivative of the smooth function
$dist^2(P,\cdot)$ has a uniform bound which does not depend on the
point $P$. Then by (5.5),(5.28) and (5.29), we have
\begin{eqnarray}
L\phi&\geqq&C_{10}\hat{C}_{18}\sum_AM(u)^{AA}-\hat{C}_{20}(1+\sum_AM(u)^{AA}),
\end{eqnarray}
where $\hat{C}_{20}$ is a constant. So by (1.2) and Proposition 3.3,
we only need to take $C_{10}$ big enough, then
 $L\phi>0$ in $\Omega^{\delta}(P)$. This means that the maximum value of function $\phi$
is attained on the boundary. Then in $\bar{\Omega}^{\delta}(P)$,
(5.32) gives
\begin{eqnarray}
(u-\psi)_C&\leqq&
C_{11}d^2+C_{10}[(1-e^{-C_8(u-\psi)})+(1-e^{-C_{9}C_8(x_1-\theta)})].
\end{eqnarray}
Both sides of the above inequality are $0$  at point $P$. Now we
obtain the uniform lower bound of term $(u-\psi)_{Cm}$. And letting
\begin{eqnarray}
\phi&=& -(u-\psi)_{C} +v,
\end{eqnarray}
we can similarly obtain the upper bound.  So for $C\neq m$, there is
a positive constant $C_{12}$ depending on $\psi$, $K$, $m$, $n$,
$\partial\Omega$ such that on $\partial \Omega_x\times S^n$
\begin{eqnarray}
|(u-\psi)_{Cm}|&\leqq& C_{12}.
\end{eqnarray}
By the choice of frame we made in the head of this section, and by
the equality on the boundary $u=\psi$, we know that for $A,B\neq m$,
on $\partial\Omega^{\delta}(P)$
\begin{eqnarray}
u_{AB}&=&\psi_{AB}-h_{AB}(u-\psi)_m,
\end{eqnarray}
where $h_{AB}$ is the second fundamental tensor along the outward
normal direction $e_m$ of $\partial\Omega$. Obviously, if one of
$A,B$ takes value in $m+1,\cdots,m+n$, then $h_{AB}=0$. Moreover, as
$\psi\in$ ABF($\tau,K$), $\psi$ only depends on $S^m$. So for
$i,j\neq m$,
\begin{eqnarray}
u_{ij}&=&\psi_{ij}-h_{ij}(u-\psi)_m,
\end{eqnarray}
and $u_{AB}=0$ for the other cases. Since $M(u)$ is positive
definite, we now only need the upper bound on $u_{mm}$. We use the
same argument as in the papers [T] and [G]. For
$P\in\partial\Omega$, define a function
\begin{eqnarray}
\lambda(P)&=&\min_{|\xi|=1,\xi\in
T_P(\partial\Omega_x)}[\nabla_{\xi\xi}u(P)-(\nabla_{\xi}u(P))^2-1],
\end{eqnarray}
where $|\cdot|$ is the standard norm of $S^m\times S^n$. Assume that
at $P_0\in
\partial\Omega$ and $\xi=e_1(P_0)\in T_{P_0}(\partial\Omega_x)$,
$\lambda$ attains its minimum value. Then $\lambda(P)\geqq
\lambda(P_0)$. By (5.37), (5.39) and $u|_{\partial\Omega}=\psi$,
\begin{eqnarray}
h_{11}(P)(u-\psi)_m(P)&\leqq&(\psi_{11}(P)-\psi_1^2(P))-(\psi_{11}(P_0)-\psi_1^2(P_0))\\
&&+h_{11}(P_0)(u-\psi)_m(P_0)\nonumber.
\end{eqnarray}
 By the compactness of $\partial\Omega_x$, there is a uniform
sufficiently small $\delta$ such that on
$\bar{\Omega}^{\delta}(P_0)$, the smooth function
$h_{11}=<\nabla_{e_1}e_1,e_m>_m$ has a uniform negative upper bound.
Now on $\bar{\Omega}^{\delta}(P_0)$, let
\begin{eqnarray}
\Psi(Q)&=&h_{11}(Q)^{-1}[(\psi_{11}(Q)-\psi_1^2(Q))-(\psi_{11}(P_0)-\psi_1^2(P_0))\\
&&+h_{11}(P_0)(u-\psi)_m(P_0)]\nonumber,
\end{eqnarray}
where $Q\in \bar{\Omega}^{\delta}(P_0)$. By (5.40), for $P\in
\bar{\Omega}^{\delta}(P_0)\cap \partial\Omega$, obviously
\begin{eqnarray}
(u-\psi)_m(P)&\geqq&\Psi(P),\ \ \text{and} \ \ (u-\psi)_m(P_0)\ \ =\
\ \Psi(P_0).
\end{eqnarray}
Then on $\bar{\Omega}^{\delta}(P_0)$, let
\begin{eqnarray}
\phi&=&\Psi-(u-\psi)_m+v.
\end{eqnarray}
By (5.42), and a similar argument as in  (5.1) to (5.5) and (5.28)
to (5.33), $\phi$ attains its maximum value at point $P_0$ for a
choice of suitable constants of $v$. Then
 $\phi_m(P_0)\geqq 0$ which implies that $u_{mm}(P_0)$ has a upper bound. So at $P_0$ by (5.38), all
eigenvalues of $M(u)$ have a upper bound. By equation (1.2), the
minimum eigenvalue of $M(u)$ at $P_0$ has a lower bound, which
implies that $\lambda (P_0)$ has a lower bound. Since $P_0$ is the
minimum value point of $\lambda$, for any boundary point $P$, and
any unit vector $\xi \in
 T_P(\partial\Omega_x)$,
$\nabla_{\xi\xi}u(P)-(\nabla_{\xi}u(P))^2-1$ has a uniform lower
bound. By the sentence after (5.38) and Definition 3.2,
\begin{eqnarray}
M(u)^{*}_{mm}&=&(u_{ij}-u_iu_j-\delta_{ij})_{(m-1)\times(m-1)},
\end{eqnarray}
where $M(u)^*_{mm}$ is the cofactor matrix of $M(u)_{mm}$, and the
right hand side is a order $m-1$ matrix with $i,j\neq m$. So on the
tangent space of $\partial\Omega_x$, diagonalizing the matrix
$M(u)^{*}_{mm}$, we find that $M(u)^*_{mm}$ has a positive uniform
lower bound. With the same argument as in [CNS1], $u_{mm}$ has a
uniform upper bound. Now by (5.36) and (5.37), there is a positive
constant $C_{13}$ depending on $\psi$, $K$, $m$, $n$,
$\partial\Omega$ such that on $\partial \Omega_x\times S^n$,
\begin{eqnarray}
|u_{AB}|&\leqq& C_{13}.
\end{eqnarray}
Now we have the $C^2$  estimate on the boundary, and combining this
with the interior $C^2$ estimate and proposition 3.3, we obtain
$C^2$ estimate. Then using Evans-Krylov theory (see[GT]), we have
the $C^{2,\alpha}$ estimate. Then differentiate equation (1.2) and
using Schauder theory, we obtain Proposition 1.2. This gives the
existence part of Theorem 1.2. For the uniqueness part, we let $g$
equals $f^{(m+n+2)/2}$ in Lemma 3.1. If $u,v$ are both solutions of
Problem (1.2), then $G(u)=G(v)$. Then similar as the argument of
(3.14),(3.15), we can use Lemma 3.1, and for the equality of
boundary-values, we have $u=v$.

\bigskip

Wang ZhiZhang

Fudan University, Institute of Mathematics Science, Shanghai,
200433, China

E-mail: youxiang163wang@163.com
\end{document}